# The system of integer functions

an efficient version of discrete mathematical analysis

József Peredy

Professor emeritus, University of Technology and Economics, Budapest, Hungary

In this short sketch the fundamental ideas of a system of discrete mathematical analysis (let's call it "the system of integer functions"), will be presented. This proposed system of integer functions is logically fully independent from the traditional mathematical analysis of the real functions, but there is a well-defined mutual correspondence between the two disciplines. The system of integer functions intends to help

- to make the transition from the present approach to the problems of the calculus to a more computer-centric one as smooth and efficient as possible[i], and
- to find a way to some kind of synthesis of the discrete and continuous[ii].

This sketch focuses on the fundamental ideas, and mainly relatively simple 2D monotonic non-decreasing functions (or such sections of functions) will be considered, but the possibilities of a wide range of generalisations will be indicated, too.

1) Integer functions

The finite sequence of $[i_k, j_k]$, (k = 0, 1, … l) ordered pairs of integers is an integer function (shortly IF), if the consecutive pairs $[i_k, j_k]$ and $[i_{k+1}, j_{k+1}]$ are neighbours in the sense, that if $j_k = j_{k+1}$ then $i_{k+1} = i_k \pm \underline{1}$ ($i^+$ or $i^-$ step), or if $i_k = i_{k+1}$ then $j_{k+1} = j_k \pm \underline{1}$ ($j^+$ or $j^-$ step).

Here the $i_k$-s and $j_k$-s are positive or negative integers (inclusive 0) and $k_l$ is a finite positive integer number. The integer pairs $[i_k, j_k]$ are the elements of the IF, the $i_k$ and $j_k$ values are the coordinates of the integer pair $[i_k, j_k]$, i. e. the variables of the IF. Consider the table below with the 25 integer pairs:

| k | 0 | 1 | 2 | 3 | 4 | 5 | 6 | 7 | 8 | 9 | 10 | 11 | 12 | 13 | 14 | 15 | 16 | 17 | 18 | 19 | 20 | 21 | 22 | 23 | 24 |
|---|---|---|---|---|---|---|---|---|---|---|----|----|----|----|----|----|----|----|----|----|----|----|----|----|----|
| i | 0 | 1 | 1 | 2 | 2 | 3 | 3 | 4 | 4 | 5 | 5 | 6 | 6 | 7 | 7 | 8 | 9 | 9 | 10 | 11 | 11 | 12 | 13 | 14 | 15 |
| j | 0 | 0 | 1 | 1 | 2 | 2 | 3 | 3 | 4 | 4 | 5 | 5 | 6 | 6 | 7 | 7 | 7 | 8 | 8 | 8 | 9 | 9 | 9 | 9 | 9 |
| lép |   | i | j | i | j | i | j | i | j | i | j | i | j | i | j | i | i | j | i | i | j | i | i | i | i |

This is really an integer function, because out of the four integers in two adjacent integer pairs either the upper two or the lower two are equal, and the other two differ by one, i. e. the integer pairs in adjacent columns are neighbours. Besides the integer pairs $[i_k, j_k]$, the table shows the serial numbers k, and the letter symbols i or j, indicating the type of the step between $[i_{k-1}, j_{k-1}]$ and $[i_k, j_k]$ as well. This integer function will be used throughout this paper as a sample.

2) The fields of differences of an integer function

In the sequence of the $i_k$ values, two consecutive values are either the same, or differ by 1. If in case of two consecutive integer pairs $[i_{k-1}, j_{k-1}]$ and $[i_k, j_k]$ the coordinate $i_k$ is by one greater or smaller than $i_{k-1}$, the integer pair $[i_k, j_k]$ will be referred to as a characteristic one: $[í_k, j_k]$. Let's consider a characteristic $[í_k, j_k]$ integer pair, and a difference class D. Find another characteristic integer pair $[í_{k'}, j_{k'}]$, where $í_{k'} = í_k + D$, and compute the characteristic difference $_D d_i = j_{k'} - j_k = k' - k - D$. The $_D d_i$ value indicates the number of j steps between $í_k$ and $í_{k'} = í_k + D$, that is it shows how many j steps belong to *D* consecutive *i* steps, starting from the location $i_k$ of the integer function. If in a *D* difference class we compute the $_Dd_i$ values for all meaningful characteristic $i_k$ – s, one after the other, it gives an overview of the relative frequencies of the *j* steps with respect to *D* consecutive *i* steps, and taking all difference classes one after the other, we get a full description of the distribution of the relative frequencies of the *i* and *j* steps resp. in the IF under investigation. The complete entirety of all $_Dd_i$ values, arranged in difference classes, and ordered to the characteristic $i_k$ values is called the field of *j* differences ordered to the characteristic $i_k$ values. The $_Dd_i$ values may be investigated in form of being "plotted" against the *j* values as well, and in view of the perfect symmetry of *i* and *j*, the field of $_Dd_j$ differences ordered either to the $j_k$ or to the $i_k$ values may be construed similarly. The table below shows the $_Dd_i$ differences ordered to the characteristic $i_k$ values of the sample IF in two difference classes. (The integer pairs belonging to characteristic $i_k$ - s are marked by green colour.)

| k | 0 | 1 | 2 | 3 | 4 | 5 | 6 | 7 | 8 | 9 | 10 | 11 | 12 | 13 | 14 | 15 | 16 | 17 | 18 | 19 | 20 | 21 | 22 | 23 | 24 |
|---|---|---|---|---|---|---|---|---|---|---|----|----|----|----|----|----|----|----|----|----|----|----|----|----|----|
| i | 0 | 1 | 1 | 2 | 2 | 3 | 3 | 4 | 4 | 5 | 5 | 6 | 6 | 7 | 7 | 8 | 9 | 9 | 10 | 11 | 11 | 12 | 13 | 14 | 15 |
| j | 0 | 0 | 1 | 1 | 2 | 2 | 3 | 3 | 4 | 4 | 5 | 5 | 6 | 6 | 7 | 7 | 7 | 8 | 8 | 8 | 9 | 9 | 9 | 9 | 9 |
| $_{D=1}d_i$ | | 1 | | 1 | | 1 | | 1 | | 1 | | 1 | | 1 | | 0 | 1 | | 0 | 1 | | 0 | 0 | 0 | |
| $_{D=3}d_i$ | | 3 | | 3 | | 3 | | 3 | | 3 | | 2 | | 2 | | 1 | 2 | | 1 | 1 | | 0 | | | |

3) Graphic representation of the IFs and the mutual correspondence with the real functions.

A visual aid for studying the IFs can be constructed if the integer pairs $[i_k, j_k]$ are represented by the unit squares $i_k \leq x < i_k + 1$ and $j_k \leq y < j_k + 1$ of a Cartesian coordinate system. Such a visual representation of the sample IF is shown on the illustration. Of course, the (x, y) Cartesian coordinate system belongs to the world real numbers and continuous functions, but the individual unit squares are discrete entities as well as the integer pairs themselves. So using such a visual aid does not interfere with the discrete character of the IFs.

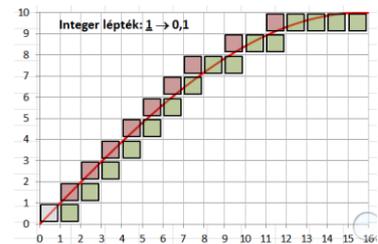

On the other hand, such a method of visual representation suggests a correct mutual correspondence between the discrete world of IFs and the world of the traditional real functions. If the unit squares are considered to be not just visual symbols but definite sets of the (x, y) points on the $\mathbb{R}^2$ plane, there are surely real functions (or finite sections of real functions) having points in each unite square belonging to a given IF and no point outside them. In such a case the real function (or a section of it) and the IF under investigation may be considered to be counterparts (a "real counterpart" or a "discrete counterpart" resp.) of each other. A real function (section) has a single, well defined IF as its discrete counterpart, while an IF has always "many" (continuum cardinality) real counterparts.

On the illustration above, in addition to the unit square symbols, a red curve appears too. This represents a finite section of a real sinus function, as a real counterpart of the sample IF. So the sample IF may be considered to be a discrete sinus function (remembering, of course, that it has many other real counterparts, too).

4) The "integer scale" and the full integer derivative

In the chapter 1 the stepwise unit modification of the *i* or *j* values is marked by double underlined $\underline{\underline{1}}$. This indicates that the unit, used for constructing an integer function may differ from the units used in "traditional" computations. It will be shown in this paper, that the system of integer functions is able to compute e. g. the values of the sine function with a given, say two digits, or four digits accuracy. In this case the unit $\underline{\underline{1}}$ of the "integer function world" will correspond with the 0,01 or 0,0001 resp. values of "real function world". The relations $\underline{\underline{1}} \to 0,01$ or $\underline{\underline{1}} \to 0,0001$ denote the "integer scales", the latter one is told to be the "finer". (Such a notation appears on the figures, representing the sample IF.)

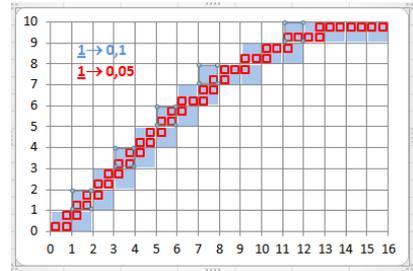

In view of this fact that the IFs may work with different scales to express the same interrelation of its two variables, one has to ensure that the results pertaining to different scales should be in accordance with each other, notably that the $_Dd_i$ differences should remain compatible while describing the same integer function with different scales. Let's consider two IFs, an original one $[i_k, j_k]$, $(k = 0, 1, \ldots l)$, and a finer, longer one $[i_k, j_k]$ $(k = 0, 1, \ldots l)$. Suppose, that the original integer scale is $\underline{\underline{1}} \to u_o$, the finer one $\underline{\underline{1}} \to u_f$, and $u_o = m.u_f$. If between the original and finer coordinates hold the inequalities $m.j_k \leq j_k < m.(j_k+1)$ (and, of course, the same for *i*), the original and the finer IFs are considered to belong together and may be used to describe similar regularities of the mutual correspondence between two variables with different scales.

In the previous chapter, the characteristic differences $_Dd_i = j_{k'} - j_k$ played a fundamental role in the definition of the fields of differences. The corresponding difference in case of the finer scale is $_Dd_i = j_{k'} - j_k$, or expressed with the original coordinate values $m.(j_{k'} - j_k) \leq j_{k'} - j_k < m.(j_{k'} - j_k +1)$. But, on the other hand, if we express the finer coordinates with the original ones first, and compute the difference after it, the result is $m.(j_{k'} - j_k -1) \leq j_{k'} - j_k < m.(j_{k'} - j_k +1)$. This means that the "normal" difference of the original scale coordinate values ($j_{k'} - j_k$) does not contain all the corresponding finer scale differences. To the "normal" difference ($j_{k'} - j_k$) the lower by 1 coordinate value ($j_{k'} - j_k -1$) has to be joined too, if both the original and the finer scales are allowed, and the contradiction free character should be preserved. So in what follows the "scaled" characteristic difference

$$_D\delta_i = \{_Dd_i, (_Dd_i-1)\}$$

will be used, that is the characteristic differences will be considered to consist out of two integers: $_Dd_i$ and $(_Dd_i-1)$.

Of course, a corresponding rule of computing differences is well known in the algebra of real intervals. The purpose of the present chapter is just to show that to arrive at this rule, it is not necessary to introduce the abstract concepts of continuity and of intervals of real numbers, it is sufficient to allow the integer numbers to represent different scales without contradiction while computing differences. So the system of integer function proposed in this paper remains logically fully independent of the concepts of continuity and real functions despite using "scaled" characteristic differences. May be, that this very simple and self-evident fact can help to establish some kind of synthesis of the discrete and the continuous approaches to the analysis of functions?

The "scaled" characteristic differences $_D\delta_i$ (or $_D\delta_j$) can be grouped into difference classes and ordered to coordinates similarly to the "normal" characteristic differences $_Dd_i$ (or $_Dd_j$). In a difference class $D$ any IF containing at least one of the two integers belonging to $_D\delta_i$ (or $_D\delta_j$) as an appropriate coordinate of one of its integer pairs, is considered to be a D class integer derivative of the IF, and a complete entirety of such integer class-derivatives for all difference classes meaningful in the (of course finite) IF is a full integer derivative. The figure below shows class derivatives of the sample IF in two difference classes D=1 and D=3:

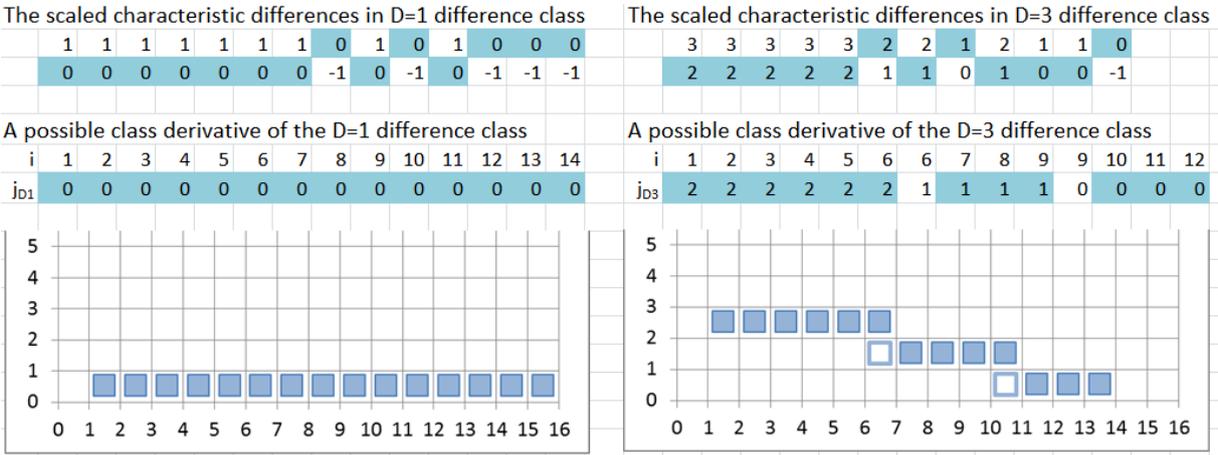

5) Algorithms generating integer functions

The core of the algorithms generating IFs is the compound conditional statement shown below:

    If RX-RY < 0 Then

        XX = XX + XXX: XY = XY + XYX: YX = YX + YXX: YY= YY+YYX

        X = X + XX: Y = Y + YX

        RX = RX + X

        i = i + 1 (i step)

    Else

        XX = XX + XXY: XY = XY + XYY: YX = YX +YXY: YY = YY +YYY

        X = X + XY: Y = Y + YY

        RY = RY + Y

        j = j + 1 (j step)

    End If.

This compound conditional statement works with two regulator registers RX and RY, and 14 work registers:

        X, Y, XX, XY, YX, YY, XXX, XXY, XYX, XYY, YXX, YXY, YYX, YYY.

All the registers take just positive or negative integer numbers (inclusive 0). The work registers have a rank, in accordance with the number of characters in their identifier.

The algorithm works step by step. In each step the algorithm first compares the actual values of the regulator registers and decides, whether the next step should be an i step or a j step, then carries on the additions prescribed for the corresponding step type, and generates the next integer pair $[i_k, j_k]$ by modifying either the previous $i_{k-1}$, or the pervious $j_{k-1}$ value by 1. A $l$ –time repetition of such steps produces the integer function required.

Of course, to generate an IF the initial values of the starting integer pair [$i$, $j$] and of all 16 registers have to be given (the initial values will be denoted by the zero index $_0$). During the process of generating an integer function, the values in many registers are modified step by step, but some of them remain constants. Such constant registers are first of all the highest rank registers (in the present case the registers: XXX, XXY, XYX, XYY, YXX, YXY, YYX and YYY registers of $3^{rd}$ rank), but some lower rank registers may remain unchanged too, if all the higher rank registers influencing their values have zero initial values. Those registers, which have non-zero initial values and remain constants, define the type of the given IF. The type of an IF is denoted by indicating the identifiers of the constant, non-zero registers in { } brackets.

6) Some discrete 2D curves.

Let's see some types of IFs. As mentioned in the previous paragraphs, the system of integer functions is logically fully independent from the calculus of real functions, but a well-defined, mutual correspondence exists between the IFs introduced in this presentation and the 2D real functions of the traditional analysis. On the basis of this correspondence, some names used for types of real functions will be applied in context with the types of IFs, too.

{ X,Y } *straight lines;* the corresponding compound conditional statement is:
    *If RX < RY Then RX = RX + X: i = i + 1 Else RY = RY + Y: j = j + 1.*

{ XX,Y } *parabolas;* the corresponding compound conditional statement is:
   If RX < RY Then X = X + XX: RX = RX + X: i = i + 1 Else RY = RY + Y: j = j + 1 .

{ XY,Y } *exponentials & logarithms;* the corresponding compound conditional statement is:
   If RX < RY Then X = X + XX: RX = RX + X: i = i + 1 Else RY = RY + Y: j = j + 1.

{ XX, YY } *ellipses & hyperbolas;* the corresponding compound conditional statement is:
   If RX < RY Then X = X + XX: RX = RX + X: i = i + 1 Else Y = Y + YY: RY = RY + Y: j = j + 1.

{ XXY, Y } *sin, cos, sh, ch and inverses;* the corresponding compound conditional statement is:
  If RX < RY Then X = X + XX: RX = RX + X: i = i + 1 Else XX = XX + XXY: Y = Y + YY: RY = RY + Y: j = j + 1 .

{ XX, YYY } *semi cubic parabolas;* the corresponding compound conditional statement is:
  If RX < RY Then X = X + XX: RX = RX + X: i = i + 1 Else YY = YY + YYY: Y = Y + YY: RY = RY + Y: j = j + 1.

The introduced type notation with { } brackets realises a full and well-arranged systematisation of "the finite world" of the discrete IFs. Under certain limitations any continuous real function may find its discrete counterpart, but in the world of the real functions there are phenomena, which are fully outside the scope of the system of IFs.

7) Possible generalisations.

As mentioned in the introduction, for the sake of lucidity the wording in this short sketch sometimes keeps in view just simple special cases. In the previous chapter, a very simple special case of the IF generating algorithm has been presented:
- The *i* or *j* steps were always positive, i. e. just monotonically non-decreasing sections were allowed. Of course, the IFs may have other sections too. The working applications harmonise

- the signs of the steps with the signs of the X and Y registers, resulting composite IFs, as seen in the figures below.
- The identifiers contained just two letter-symbols X and Y, sufficient for drawing 2D curves. It is possible to use X, Y and Z to represent curves or surfaces in 3D, or work in the space-time using X, Y, Z and T. The basic structure of the IF algorithm remains unchanged even in these more general cases: there are three or four kinds of steps, say *i*, *j*, *z*, or even *t*, and each work register is in function in those steps, which correspond with the last letter-symbol of its identifier. Being in function means that it has to be added to the (by one lower rank) register, won by dropping the last letter of the identifier. Anyway, it has to be mentioned that to find the appropriate sequence of the consecutive *i*, *j*, *z* (or even *t*) types of the steps, needs somewhat modified method in the higher dimensions compared with the simple investigation of the RX-RY < 0 inequality, as shown in chapter 5.
- In chapter 5. the highest rank among the work registers was 3. As the previous samples show, even this rank allows representing a fairly wide variety of different types of IFs, but the rank may be raised, generally may be any finite positive integer.

See below, as examples of non-monotonic discrete 2D curves, a semi-cubic egg curve and a section of a sinusoidal curve, drawn by a software application of the IF generating algorithm:

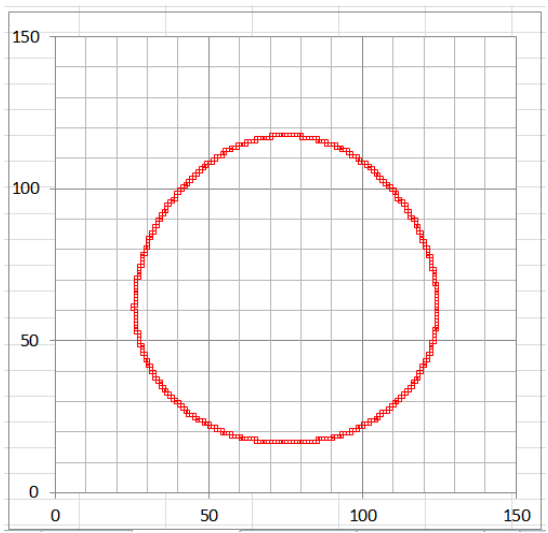

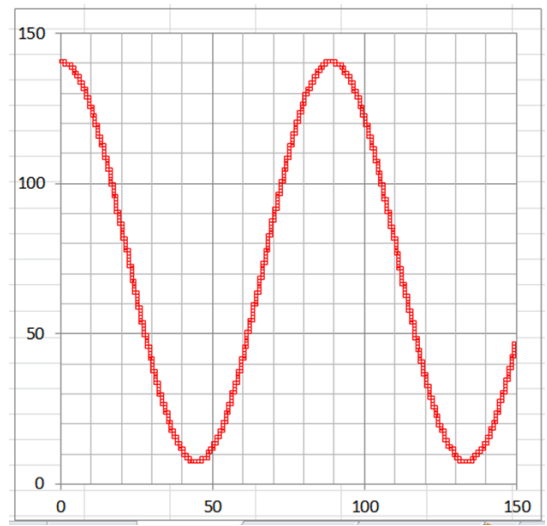

Type: *{ XX, YYY }*. The actual initial values:

| i= | 25 | j= | 60 | R= | 0 | k= | 2000 |
|---|---|---|---|---|---|---|---|
| X | 500000 | Y | 10 | | | | |
| XX | -10000 | XY | 0 | YX | 0 | YY | 10000 |
| XXX | 0 | XXY | 0 | XYX | 0 | XYY | 0 |
| YXX | 0 | YXY | 0 | YYX | 0 | YYY | -125 |

Type: *{ XXY, Y }* . The actual initial values:

| i= | 0 | j= | 140 | R= | 500 | 1 | 2000 |
|---|---|---|---|---|---|---|---|
| X | 0 | Y | 600 | | | | |
| XX | -200 | XY | 0 | YX | 0 | YY | 0 |
| XXX | 0 | XXY | -3 | XYX | 0 | XYY | 0 |
| YXX | 0 | YXY | 0 | YYX | 0 | YYY | 0 |

Of course, the computer applications of the IF algorithm might work in full accordance with the resolution defined by hardware; i. e. could draw as "smooth" curves as possible on the given display. The above presented figures intend to emphasize the discrete character of the IFs.

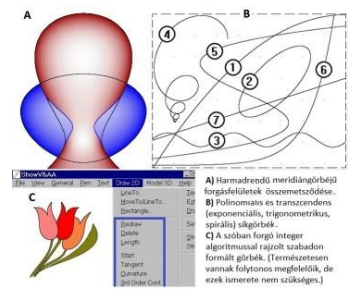

There are computer applications of the system of integer functions for 3D geometry, and for drawing freeform curves, too (see right).

8) Integer functions and the corresponding parallel sequences of the regulator register values.

The algorithms generating integer functions in each step of their functioning produce either a new RX value (in i steps), or a new RY value (in j steps). These RX and RY values form two separate sequences of integer numbers

$$RX_0, RX_1, \ldots RX_i, \ldots RX_{il} \text{ and } RY_0, RY_1, \ldots RY_j, \ldots R_{jl}.$$

In this introductory paper monotonic non-decreasing RX and RY sequences will be considered only (that is such IFs or such sections of IFs where this restriction is valid). The generating algorithm produces the $RX_i$ and $RY_j$ values in the order of the consecutive steps, that is in a single combined sequence, where they have a second index $k$ ($RX_{i,k}$ or $RY_{j,k}$) specifying the serial number of the step in which they were generated. The combined sequence in its entirety is no more monotonic non-decreasing, but if we consider just the characteristic integer pairs and the integer pairs immediately preceding them, this "incomplete" sequence preserves the monotonic non-decreasing character. Such pairs of regulator value sequences (R sequences) are of use for studying "theoretically" different properties of the IFs.

9) Some elementary problems of mechanics and the system of the integer functions.

"The main external source of mathematical problems is science, in particular physics."[iii] Let's see, how the system of integer function works in case of some elementary problems of mechanics. Such problems were studied by Newton and Leibnitz, while doing the first steps towards the modern foundation of mathematical analysis.

9.1 <u>The Newton's First Law</u>:
> "every body perseveres in its state of being at rest or of moving uniformly straight
> forward except insofar as it is compelled to change its state by forces impressed".[iv]

Such "uniformly straight forward" motions can be described by simple integer functions. If we have a measurement resulting that during the time of $i_l$ time units, the covered distance is $j_l$ linear measure units, the corresponding IF has to start with $[i_0, j_0] = [0, 0]$, and to end with $[i_l, j_l]$. This means, that the corresponding IF contains $i_l$ times i steps and $j_l$ times j steps, distributed as "uniformly", as possible to describe a uniform motion.

An IF of type $\{X,Y\}$ can realize this, by choosing the initial values $XR_0 = YR_0 = 0$, $X = j_l$ and $Y = i_l$. In this case the R values are

$$RX_i = i.X \text{ and } RY_j = j.Y,$$

and the characteristic difference $_Dd_i$ belonging to a characteristic $i$ and a given difference class $D$ is the greatest value fulfilling the inequality

$$i.X + D.X - j.Y - _Dd_i.Y \geq 0,$$

because the $_Dd_i$ is the difference of $j$ coordinates of two *characteristic* elements. The same reason defines, that $0 \leq (j.Y - i.X) < Y$. Consequently, the characteristic differences are either $_Dd_i = Int(D.X/Y)$, or $_Dd_i = Int(D.X/Y) + 1$, and the "scaled" characteristic differences $_D\delta_i = \{ _Dd_i, (_Dd_i-1)\}$ contain the $Int(D.X/Y)$ value in both possible cases. So the constant integer functions

$$[i_k, Int(D.X/Y)], (k = 0, 1, \ldots l-D)$$

are integer derivatives in the corresponding *D* classes, and the entirety of all these *(l-1)* items of class derivatives constitutes a full integer derivative of the IF of type ⸤ *X, Y* ⸥. That is in all difference classes exists a constant integer derivative, and this means that the difference between two $_Dd_i$ values belonging to any i may be just 1, if any, that is the i steps and j steps are distributed really as uniformly as possible.

The $_Dd_i$ characteristic differences show, how many linear measure units are done during *D* time units of the uniform motion, i. e. if Y is constant, the greater is the X value the greater is the velocity.

9.2 The free fall, as explained in the 3$^{rd}$ edition of the Principia. It is interesting to observe, that in this text one may find the expression "individual equal particles of time", i. e. some traces of the discrete attitude are present:

> "When a body falls, uniform gravity, by acting equally in individual equal particles of time, impresses equal forces upon that body and generates equal velocities; and in the total time … it generates a total velocity proportional to the time."

The system of integer functions lends itself to follow this train of thought of Newton. In case of free fall of a body one investigates the correlation between the elapsed time and the covered distance from the starting point of time and from the starting position. The i steps represent the "individual equal particles of time" and the "generated" additional "equal velocities" can be represented by giving a non-zero initial value to the XX register, which will be added in each i step to the velocity, represented by X (in case of a constant Y). That is the free fall can be described by integer functions of type ⸤*XX, Y*⸥, where XX represents the change of the velocity X in a time unit, that is the acceleration.

9.3 Simple harmonic motion. As a third elementary problem let's see the simple harmonic motion as 
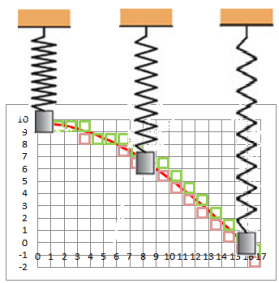
it appears when a linearly elastic spring hinders the free fall. In free fall "equal forces generate equal velocities". In the presence of an elastic spring the forces depend on the actual distance of the "body" from its starting position in such a manner that equal changes of the position (j steps) generate equal changes of the force and consequently equal changes of the acceleration XX. In each j step the XX register has to be modified by a constant XXY value, that is the initial value of the work register XXY is non-zero, and the simple harmonic motion can be described by integer functions of type ⸤ *XXY, Y* ⸥.

At least in the case of the above three elementary problems the system of IFs proved to be usable in mathematical description of physical phenomena.

10) The "discrete π".

In the traditional mathematical analysis the simple harmonic motion is described by sinusoidal real functions. In these functions the value π has an important role. If the IFs of type ⸤ *XXY, Y* ⸥ are really discrete counterparts of the real sinusoidal functions (among others), the π value must appear in some form even in the discrete world.

Let's consider the $y = sin(x)$ real function. It passes through the origin, having here a tangent with $45^0$ degree slope and a zero curvature (as a point of inflexion). The next characteristic point is its maximum point, with $x = π/2$ and $y = 1$ coordinates and a horizontal tangent.

Let's try to construct now the same situation in the world of IFs. The initial integer pair should be $[i_0=0, j_0=0]$, the initial values of the work registers X and Y should be equal positive numbers $X_0 = Y_0$ (to ensure the $45^0$ degree tangent), and the XX work register should change its sign to negative in the first step (XX = 0, XXY =-1). Starting the ɫ XXY, Y ɫ IF algorithm with such initial values let it run until that step, in which the work register X changes its sign to negative. Let denote $[i_{π/2}, j_{π/2}]$ the integer pair belonging to this step. Remembering, that the corresponding coordinate in the "real world" is $y = 1$, the integer scale is $\underline{\underline{1}} \rightarrow (1/j_{π/2})$. On the basis of this scale one can go over into the world of the real, continuous functions, and using the unit squares belonging to the integer pairs [0, 0] and $[i_{π/2}, j_{π/2}]$ resp. can compute real upper and lower limits for the real value of $π/2$. Of course, the integer pair $[i_{π/2}, j_{π/2}]$ is exact, but the limits expressed as finite decimal fractions are approximations.

In the previous paragraph appeared the requirement $X_0 = Y_0$, but no mention has been made of the common actual value. As a matter of fact, any positive integer value gives correct results for the limits, but the greater is the $X_0 = Y_0$ initial value, the narrower (i. e. the more accurate) will be the limit. In the table below the results of a few experiments are listed. The computations run on a commercial laptop, with a simple Excel-VBA program (see APPENDIX).

| X = Y | $i_{π/2}$ | $j_{π/2}$ | alsó korlát | felső korlát | idő(sec) |
|---|---|---|---|---|---|
| 100 | 15 | 8 | 1,4 | 1,777778 | - |
| 10000 | 157 | 99 | 1,56 | 1,59596 | - |
| 10000000 | 4967 | 3161 | 1,570524 | 1,571655 | 0,01 |
| 3,7683E+11 | 304924 | 194120 | 1,570788 | 1,570807 | 0,05 |
| 1E+12 | 1570796 | 999999 | 1,570795 | 1,570799 | 0,24 |
| 1E+14 | 15707963 | 9999999 | 1,570796 | 1,570797 | 2,21 |

"Traditional" $π/2 = 1,570\ 796\ 326\ 794…$  **lover bound   upper bound   runtime**

The last version has given 6 decimal digits accuracy. During the 2 sec. runtime more than 15 million $sin(x)$ values and ≈10 million $arcsine(y)$ values have been determined with this accuracy to arrive at the required $[i_{π/2}, j_{π/2}]$ result.

11) The system of IFs, as a new approach to the world of mathematical functions.

In the system of IFs some problems, which are well known in the "traditional" real analysis simply disappear. An IF expresses the mutual correspondence between two variables, having exactly the same role, so there are no special inverse functions: any IF is in the same time its own inverse. Consequently the whole problematic of finding the roots of an $f(x, y) = 0$ equation reduces just to count normally that integer pair(s) in which the specified coordinate has zero value, e. g. that is true for the quantic equations, too (their integer counterparts are IFs of type ɫXXXXX, Yɫ. The question, whether a function can be expressed in closed form, or can't, doesn't exist either, instead of that one may examine the highest rank of the work registers in the ɫ ɫ "struck-through square brackets" type notation, the lower this rank, the simpler is the rule of the mutual correspondence of the two variables expressed by the given IF.

The system of IFs is not just a practical method of computing, or of evading awkward problems of the real analysis, but it has certain theoretical consequences of its own, too. Perhaps the most fundamental is that the derivatives in the real analysis are (limits of) *quotients*, while the integer class derivatives are *differences*. Consequently, the absolute values of the class derivatives, showing the relative frequencies of the j (or i) steps, are (especially in the lover difference classes of long IFs) definitely less than the original j (or i) coordinates. The class derivatives themselves are IFs, so they have their class derivatives, too (say these are the second order class derivatives of the original IF). Repeating this process a few times (depending on the rank), the absolute values of the higher order, but lover difference class derivatives are more and more decreasing and finally don't contain efficient information any more. This remarkable property of the system of IFs does not appear among the higher order derivatives of the real functions, and may seem at the first sight rather strange. However, it is in close relation with the fact, that in case of discrete counterparts of the usual smooth real curves, the lower rank registers have greater initial values than the higher rank ones (see the two discrete 2D curve above), and all these shed light on the fact, that the different factors, expressed by the initial values of the different work registers influence the magnitude of the coordinates of the IFs in different degree.

12) Closing remarks.

The system of IFs is very rich, and even in its present early initial stage came up many interesting features, which couldn't be mentioned in this short sketch. Perhaps the most remarkable is that no discrete counterpart has been presented for the "first four rule of arithmetic". As a matter of fact, all this four rules of arithmetic are included in the IFs of type *ɟX, Yɟ*, but no doubt, if the system of IFs is to be used in everyday computations, a special "discrete multiplication-division" is needed. This can be deduced from the step structure of the *ɟX, Yɟ* type IFs. The result is in itself interesting: the "discrete multiplication-division" is just a generalised form of the place value number systems. It seems that the system of IFs may open new vistas of research, of software development and even of new ideas in the processor architecture.

---

[i] https://quotefancy.com/quote/1755684 /Doron-Zeilberger
[ii] László Lovász: Discrete or Continuous: Two sides of the same?
[iii] László Lovász: Discrete or Continuous: Two sides of the same?
[iv] Michael Spivak: Elementary Mechanics from a Mathematician's Viewpoint

APPENDIX

```
Sub picomp()
Dim R As LongLong, X As LongLong, XX As LongLong, XXY As LongLong, Y As LongLong
Dim i As Long, j As Long
Dim Start, Finish, TotalTime

With Worksheets("Csemege PI()")
   Start = Timer

   XXY = -1: XX = -1: X = .Cells(17, 3).Value: Y = X
   R = 0: i = 0: j = 0

   Do
      If R > 0 Then
           XX = XX + XXY: R = R - Y: j = j + 1 'j lépés
           Else
           X = X + XX: R = R + X: i = i + 1 'i lépés
      End If
   Loop While X > 0

   .Cells(19, 2).Value = i: .Cells(19, 5).Value = j

   Finish = Timer       ' Set end time.
   TotalTime = Finish - Start    ' Calculate total time.
   MsgBox "Runtime " & TotalTime & " seconds"

End With
End Sub
```